\documentclass[12pt]{article}
\usepackage{amssymb,amsthm,amsmath,amsfonts,latexsym,tikz,hyperref}
\usepackage[hmargin=1in,vmargin=1in]{geometry}

\newtheorem{thm}{Theorem}[section]
\newtheorem{prop}[thm]{Proposition}
\newtheorem{cor}[thm]{Corollary}
\newtheorem{lem}[thm]{Lemma}
\newtheorem{conj}[thm]{Conjecture}
\newtheorem{exa}[thm]{Example}

\DeclareMathOperator{\cen}{cen}
\DeclareMathOperator{\mi}{mi}
\DeclareMathOperator{\GK}{GK}
\DeclareMathOperator{\GKh}{\widehat{GK}}
\DeclareMathOperator{\core}{core}

\newcommand{\ben}{\begin{enumerate}}
\newcommand{\een}{\end{enumerate}}
\newcommand{\ble}{\begin{lem}}
\newcommand{\ele}{\end{lem}}
\newcommand{\bth}{\begin{thm}}
\renewcommand{\eth}{\end{thm}}
\newcommand{\bpr}{\begin{prop}}
\newcommand{\epr}{\end{prop}}
\newcommand{\bco}{\begin{cor}}
\newcommand{\eco}{\end{cor}}
\newcommand{\bcon}{\begin{conj}}
\newcommand{\econ}{\end{conj}}
\newcommand{\bde}{\begin{defn}}
\newcommand{\ede}{\end{defn}}
\newcommand{\bex}{\begin{exa}}
\newcommand{\eex}{\end{exa}}
\newcommand{\barr}{\begin{array}}
\newcommand{\earr}{\end{array}}
\newcommand{\btab}{\begin{tabular}}
\newcommand{\etab}{\end{tabular}}
\newcommand{\beq}{\begin{equation}}
\newcommand{\eeq}{\end{equation}}
\newcommand{\bea}{\begin{eqnarray*}}
\newcommand{\eea}{\end{eqnarray*}}
\newcommand{\bal}{\begin{align*}}
\newcommand{\bce}{\begin{center}}
\newcommand{\ece}{\end{center}}
\newcommand{\bpi}{\begin{picture}}
\newcommand{\epi}{\end{picture}}
\newcommand{\bpp}{\begin{picture}}
\newcommand{\epp}{\end{picture}}
\newcommand{\bfi}{\begin{figure} \begin{center}}
\newcommand{\efi}{\end{center} \end{figure}}
\newcommand{\bprf}{\begin{proof}}
\newcommand{\eprf}{\end{proof}\medskip}
\newcommand{\capt}{\caption}
\newcommand{\bsl}{\begin{slide}{}}
\newcommand{\esl}{\end{slide}}
\newcommand{\bfr}{\begin{frame}}
\newcommand{\efr}{\end{frame}}

\newcommand{\hqed}{\hfill \qed}

\newcommand{\eqqed}[1]{\vspace{5pt}$\rule{1ex}{0ex}\hfill{\dil#1}\hfill\qed$}

\newcommand{\hso}[1]{\hspace{-1pt}}

\newcommand{\qmq}[1]{\quad\mbox{#1}\quad}

\newcommand{\emp}{\emptyset}

\newcommand{\sbe}{\subseteq}

\newcommand{\iso}{\cong}

\newcommand{\oh}{\hat{1}}

\newcommand{\lt}{\lhd}
\newcommand{\gt}{\rhd}
\newcommand{\lte}{\unlhd}

\newcommand{\case}[4]{\left\{\barr{ll}#1&\mbox{#2}\\#3&\mbox{#4}\earr\right.}
\newcommand{\fl}[1]{\lfloor #1 \rfloor}
\newcommand{\ce}[1]{\lceil #1 \rceil}

\def\<{\langle}
\def\>{\rangle}

\newcommand{\ree}[1]{(\ref{#1})}

\newcommand{\ra}{\rightarrow}

\newcommand{\al}{\alpha}
\newcommand{\be}{\beta}
\newcommand{\ga}{\gamma}

\newcommand{\ka}{\kappa}

\DeclareMathOperator{\rk}{rk}

\newcommand{\dil}{\displaystyle}

\begin{document}
\pagestyle{plain}

\title{On a rank-unimodality conjecture of  Morier-Genoud and Ovsienko
}
\author{Thomas McConville\\[-5pt]
\small Department of Mathematics, Kennesaw State University\\[-5pt]
\small Kennesaw, GA 30144, USA, {\tt tmcconvi@kennesaw.edu}\\
Bruce E. Sagan\\[-5pt]
\small Department of Mathematics, Michigan State University,\\[-5pt]
\small East Lansing, MI 48824, USA, {\tt sagan@math.msu.edu}\\
Clifford Smyth\\[-5pt]
\small Department of Mathematics and Statistics, University of North Carolina\\[-5pt]
\small Greensboro, NC 27402, USA, {\tt cdsmyth@uncg.edu}
}

\date{\today\\[10pt]
	\begin{flushleft}
	                                       \\[5pt]
	\end{flushleft}}

\maketitle

\begin{abstract}
Let $\al=(a,b,\ldots)$ be a composition.  Consider the associated poset $F(\al)$, called a fence, whose covering relations are
$$
x_1\lt x_2 \lt \ldots\lt x_{a+1}\gt x_{a+2}\gt \ldots\gt x_{a+b+1}\lt x_{a+b+2}\lt \ldots\ .
$$
We study the associated distributive lattice $L(\al)$ consisting of all lower order ideals of $F(\al)$.  
These lattices are important in the theory of cluster algebras and their rank generating functions can be used to define $q$-analogues of rational numbers.
In particular, we make progress on a recent conjecture of Morier-Genoud and Ovsienko that $L(\al)$ is rank unimodal.  We show that if one of the parts of $\al$ is greater than the sum of the others, then the conjecture is true.  We conjecture that $L(\al)$ enjoys the stronger properties of having a nested chain decomposition and having a rank sequence which is either top or bottom interlacing, the latter being a recently defined property of sequences.  We verify that these properties hold for compositions with at most three parts and for what we call $d$-divided posets, generalizing work of Claussen and simplifying a construction of Gansner.
\end{abstract}

\section{Basic definitions and background}

We will be studying the conjectured rank unimodality of certain distributive lattices.  We begin by defining the posets from which they arise.
Our terminology for partially ordered sets and other structures will follow Sagan's combinatorics text~\cite{sag:aoc}.
Let $\al=(\al_1,\al_2,\ldots,\al_s)$ be a {\em composition of $n-1$}, that is, a sequence of positive integers summing to $n-1$.  To simplify notation we will sometimes write $\al=(a,b,c,\ldots)$.
For each $\al$ we have a corresponding {\em fence poset}, $F=F(\al)$, with elements $x_1,x_2,\ldots,x_n$ and covering relations
\beq
\label{F(al)}
x_1\lt x_2 \lt \ldots\lt x_{a+1}\gt x_{a+2}\gt \ldots\gt x_{a+b+1}\lt x_{a+b+2}\lt\ldots\lt x_{a+b+c+1}\gt x_{a+b+c+2} \gt\ldots
\eeq
where $\lte$ is the order relation in $F$.
The Hasse diagram of the fence $F(2,3,1)$ is shown in Figure~\ref{F(2,3,1)}.
We will call the maximal chains of $F$ {\em segments} so that the $i$th part of $\al$ is equal to the length of the $i$th segment of $F$.  Because of this convention, the sum of the parts of $\al$ is one less than $\#F$, the cardinalty of $F$.  We will also denote cardinality by using the absolute value symbol.  

\bfi
\begin{tikzpicture}
\fill(0,1) circle(.1);
\fill(1,2) circle(.1);
\fill(2,3) circle(.1);
\fill(3,2) circle(.1);
\fill(4,1) circle(.1);
\fill(5,0) circle(.1);
\fill(6,1) circle(.1);
\draw (0,1)--(2,3)--(5,0)--(6,1);
\draw (0,.5) node{$x_1$};
\draw (1,1.5) node{$x_2$};
\draw (2,2.5) node{$x_3$};
\draw (3,1.5) node{$x_4$};
\draw (4,.5) node{$x_5$};
\draw (5,-.5) node{$x_6$};
\draw (6,.5) node{$x_7$};
\end{tikzpicture}
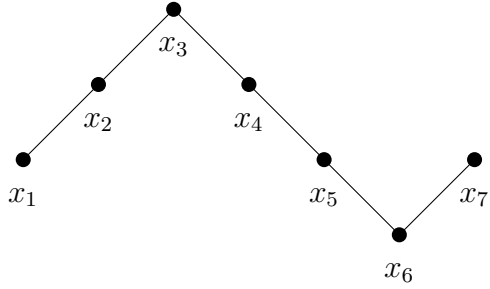
\capt{The fence $F(2,3,1)$ \label{F(2,3,1)}}
\efi

Given any poset $P$, its set of (lower) order ideals forms a distributive lattice $L(P)$.  We will shorten $L(F(\al))$ to $L(\al)$ and use similar abbreviations with other notation.  The complement of an order ideal of $P$ is an order ideal of $P^*$, the poset dual of $P$.  And 
if $\al=(\al_1,\ldots,\al_s)$ has an odd number of segments then $F(\al^r)\iso F(\al)^*$ where where $\al^r=(\al_s,\ldots,\al_2,\al_1)$ is the reversal of $\al$ and  $\iso$ is poset isomorphism.  Combining these two observations and translating to the corresponing lattices we have the following result which we record for future use.
\ble
\label{al^r}
For any $\al=(\al_1,\ldots,\al_s)$ with $s$ odd we have

\eqqed{
L(\al) \iso L(\al^r)^*.
}
\ele

The lattices $L(\al)$ will be our principal objects of study.  They are important objects in the theory of cluster algebras.  In particular, one can view $F(\al)$ as a quiver formed from the Dynkin diagram of type $A$ by replacing each cover $x\lt y$ with an arrow from $x$ to $y$.  
Then $L(\al)$ can be used to compute a mutation in a corresponding cluster algebra on a surface.  In fact, there are (at least) six different descriptions of $L(\al)$ or its dual which are useful for this computation.  These are in terms of perfect matchings on
snake graphs~\cite{pro:cfp}, 
perfect matchings of angles~\cite{yur:cef,yur:cce},
$T$-paths~\cite{sch:cef,sch:caa,ST:caa},
lattice paths on snake graphs~\cite{pro:cfp},
lattice paths of angles~\cite{cla:epp},
or $S$-paths~\cite{cla:epp}.

In order to introduce the conjecture on which we will focus, we need some definitions related to sequences and their generating functions.  We say that a sequence of nonnegative real numbers $a_0,a_1,\ldots,a_n$ is {\em unimodal} if there is some index $m$ such that
$$
a_0\le a_1\le \ldots \le a_m\ge a_{m+1}\ge \ldots \ge a_n.
$$
Unimodal sequences arise frequently in combinatorics, algebra, and geometry; see the 
survey articles of Stanley~\cite{sta:lus}, Brenti~\cite{bre:lus} or Br\"and\'en~\cite{bra:ulr}.  We will say that the generating function $f(q)=\sum_k a_k q^k$ has a property such as unimodality if its coefficient sequence does.

Now suppose that $P$ is a finite  poset.  We call $P$ {\em ranked} if, for all $x\in P$, the length of any saturated chain from a minimal element to $x$ is invariant.  This length is called the {\em rank of $x$} and denoted
$\rk x$.  We also define the {\em rank of $P$}, $\rk P$, to be the maximum of $\rk x$ over all $x\in P$. The {\em $k$th rank} of $P$ is the set
$$
R_k(P) = \{x\in P \mid \rk x=k\}
$$
and we let $r_k(P)=\#R_k(P)$.  Any finite distributive lattice is ranked by the cardinality of each element viewed as an order ideal of the corresponding poset of join irreducibles.  We say that $P$ is {\em rank unimodal} if the sequence $r_0(P), r_1(P),\ldots,r_n(P)$ is unimodal where $n=\rk P$.  We will similarly prepend ``rank" to other properties of sequences when applied to the rank sequence of a poset.  Our main object of study is the following conjecture of Morier-Genoud and Ovsienko.
\bcon[\cite{MGO:qdr}]
\label{MGO}
For any $\al$, the lattice $L(\al)$ is rank unimodal.
\econ

We note that Morier-Genoud and Ovsienko used the rank generating functions for the $L(\al)$ to define $q$-analogues for rational numbers.  Interestingly, special cases of this conjecture had already been proven even before it was stated because the problem is so natural in its own right.  Gansner~\cite{gan:loi} proved Conjecture~\ref{MGO} for certain dual fences which we call $d$-divided and which will be precisely defined in Section~\ref{ddf}.  Munarini and Zagaglia~\cite{MZ:rpl} gave a different proof of the conjecture for $2$-divided fences which are those with 
$\al=(1,1,\ldots,1)$.  Since the conjecture was posed, Claussen~\cite{cla:epp} has shown that it is true for all fences with at most four segments.  One of our main results is that Conjecture~\ref{MGO} holds if one of the segments is sufficiently long.
\bth
\label{long}
Suppose $\al=(\al_1,\ldots,\al_s)$ and there is an index $t$ such that 
$$
\al_t > \sum_{i\neq t} \al_i.
$$
Then $L(\al)$ is rank unimodal.
\eth

We will also be interested in various strengthenings of Conjecture~\ref{MGO}.  To state them, we will need to define other properties of sequences and posets.  Say that the sequence $a_0,a_1,\ldots,a_n$ is {\em symmetric} if $a_k = a_{n-k}$ for $k<n/2$.  Symmetric unimodal sequences are common, for example a row of Pascal's triangle or the coefficients of a $q$-binomial coefficient.  Even when one does not have symmetry, there may be some relation between $a_k$ and $a_{n-k}$.  Call the sequence {\em top heavy} (respecively, {\em bottom heavy}) if $a_k\le a_{n-k}$ (respectively, $a_k\ge a_{n-k}$) for $k<n/2$.   As an illustration, a special case of a result of Bj\"orner and Ekedahl~\cite{BE:sbi} states that the rank sequence for Bruhat order on a crystallographic Coxeter group is top heavy.
More recently, a new property of sequences has been identified which implies both uimodality and heaviness.  Call the sequence {\em top interlacing} if
\beq
\label{ti}
a_0\le a_n \le a_1\le a_{n-1} \le \ldots\le a_{\ce{n/2}}
\eeq
where $\ce{\cdot}$ is the ceiling function.  Similarly, the sequence is {\em bottom interlacing} if
$$
a_n\le a_0 \le a_{n-1} \le a_1 \le \ldots \le a_{\fl{n/2}}
$$
with $\fl{\cdot}$ being the floor function.  See~\cite{ath:esl,ath:gcg,BJM:hpz,BS:sdr,SV:uqi,sol:sns} 
for research related to this concept.  We note that in the literature~\ree{ti} has been called ``alternately increasing."  However, we prefer our terminology both because ``alternating" usually refers to a sequence satisfying $a_0<a_1>a_2<a_3>\ldots$, and since~\ree{ti} implies that the first half of the sequence and the reverse of the second half interlace in the usual sense of the term.  We propose the following strengthening of Conjecture~\ref{MGO}.  In it, we refer to the rank sequence 
\beq
\label{r(al)}
r(\al) = (r_0(\al),r_1(\al),\ldots,r_n(\al))
\eeq
where $r_k(\al)=r_k(L(\al))$ and $n=\# F(\al)$.  
\bcon
\label{heavy}
Suppose $\al=(\al_1,\ldots,\al_s)$.
\ben
\item[(a)]  If $s=1$ then $r(\al) = (1,1,\ldots,1)$ is symmetric.
\item[(b)]  If $s$ is even, then $r(\al)$ is bottom interlacing.
\item[(c)]  Suppose $s\ge3$ is odd and let $\al'=(\al_2,\ldots,\al_{s-1})$.
	\ben
	\item[(i)] If $\al_1>\al_s$ then $r(\al)$ is bottom interlacing.
	\item[(ii)] If $\al_1<\al_s$ then $r(\al)$ is top interlacing.
	\item[(iii)] If $\al_1=\al_s$ then $r(\al)$ is symmetric, bottom interlacing, or top interlacing depending on whether 
	$r(\al')$ is symmetric, top interlacing, or bottom interlacing, respectively.
	\een
\een
\econ

Statement (a) in this conjecture is trivial, but is needed as a base case.  
We have verified this conjecture by computer for up to $5$ segments of lengths at most $10$, and for $6$ segments having lengths at most $5$. 
We have been able to prove the conjecture for various fences, including those with at most three segments and the $d$-divided posets, by showing that the corresponding lattices satisfy an even stronger condition which we now describe.

Let $P$ be a ranked poset with $\rk P = n$.  Also, let $C$ be a saturated $x$--$y$ chain in $P$.  The {\em center} and {\em interval} of $C$ are
the rational number
$$
\cen C =\frac{\rk x + \rk y}{2}
$$
and interval of integers
$$
[C]=[\rk x,\rk y],
$$
repectively
A {\em chain decomposition} or {\em CD} of $P$ is a set partition of $P$ into saturated chains.  In a {\em symmetric chain decomposition} or {\em SCD}, every chain $C$ in the partition must satisfy $\cen C = n/2$.  Equivalently, if $C$ is an $x$--$y$ chain of the partition then  $\rk y = n-\rk x$.  If $P$ admits an SCD then its rank sequence is symmetric and unimodal.  In fact, $P$ even enjoys the strong Sperner property which says that, for all $k\ge1$,  the maximum number of elements in a subposet whose longest chain has at most $k$ elements is just the sum of the $k$ largest ranks.  See the survey article of Greene and Kleitman~\cite{GK:ptt} for more information about chain decompositions and the Sperner property.  To deal with the case when the rank sequence is not symmetric,  consider a 
{\em nested chain composition}, NCD, which is a CD where any two of its chains $C,D$ satisfy either $[C]\sbe[D]$ or $[D]\sbe[C]$.  If $P$ admits an NCD then it is rank unimodal and still has the strong Sperner property.  We will be particularly concerned with a special type of NCD.  Call a CD {\em top centered} if every chain $C$ in the partition satisfies $\cen C =n/2$ or $\cen C = (n+1)/2$.  It follows easily  that this is an NCD and the rank sequence of $P$ is top interlacing.  
Similarly, a CD is {\em bottom centered} if its chains satisfy 
$\cen C = n/2$ or $(n-1)/2$.  Again, this is an NCD and the rank sequence is now bottom interlacing.  
Note also that if a poset has an NCD and its rank sequence is top or bottom interlacing then the NCD must be top or bottom centered, respectively.
This can be proven inductively using the observation that an NCD must contain a chain from a minimum rank element to a maximum rank element.
This leads to the strongest of our conjectures so far.
\bcon
\label{centered}
For any $\al$, the lattice $L(\al)$ admits a CD which is either symmetric, top centered, or bottom centered consistent with Conjecture~\ref{heavy}.
\econ

 We prove a number of special cases of this conjecture in the sequel.  In particular, when the fence has at most three segments we have the following refinement of Claussen's result on the rank unimodality of the $L(\al)$.
\bth
\label{<=3}
If $\al$ has at most three parts then Conjecture~\ref{centered} is true.
\eth

The rest of this paper is structured as follows.  In the next section we will prove Theorem~\ref{long} in the case that the long segment is the first or the last.  We will do this using a recursion which will also permit us to replace the strict inequality with a weak one for these particular segments.  In Section~\ref{lms} we will complete the proof of Theorem~\ref{long}.  We will also describe an inductive procedure for proving Conjecture~\ref{centered} when a long segment exists.  The following section will be devoted to proving Theorem~\ref{<=3}.     Section~\ref{ddf} contains a proof of Conjecture~\ref{centered} for $d$-divided posets.  It is modeled on, but simpler than, Gansner's proof.  We end with a section outlining two more possible approaches to these conjectures.

\section{Long initial or final segments}
\label{lifs}

In this section we will prove a stronger version of Theorem~\ref{long} where the long segment is either the first or the last.  This will be based on a recursion for the rank generating function
$$
r(q;\al) = \sum_{k\ge0} r_k(\al) q^k
$$
where the $r_k(\al)$ are given by~\ree{r(al)}.  The method of proof involves considering the ideals of $F(\al)$ which do or do not contain an element $x$ which we will call {\em toggling on $x$}.  Also, for the recursions to make sense, we must permit compositions
$(\al_1,\ldots,\al_s)$ where $\al_s=0$.  But in this case we just define
$$
F(\al_1,\ldots,\al_{s-1},0)=F(\al_1,\ldots,\al_{s-1}).
$$
\ble
\label{al_sRec}
Let $\al=(\al_1,\al_2,\ldots,\al_s)$.  Then for $s$ odd
$$
r(q;\al)=r(q;\al_1,\ldots,\al_{s-1},\al_s-1)  + q^{\al_s+1}\cdot r(q;\al_1,\ldots,\al_{s-2},\al_{s-1}-1)
$$
and for $s$ even
$$
r(q;\al) = r(q;\al_1,\ldots,\al_{s-2},\al_{s-1}-1) +  q\cdot r(q;\al_1,\ldots,\al_{s-1},\al_s-1).
$$
\ele
\bprf
Induct on $s$.  We will just do the induction step for $s$ odd as the case when $s$ is even is similar.  Let $x=x_n$ be the furthest right element of $F(\al)$ in~\ree{F(al)}.  Then a given lower order ideal $I$ of $F(\al)$ contains $x$ or does not contain $x$.  In the case where $I$ does not contain $x$, the fact that $s$  is odd implies that these ideals are in bijection with the ideals of $F(\al_1,\ldots,\al_{s-1},\al_s-1)$.  So the contributions of such order ideals to $r(q;\al)$ give the first summand in the recursion.  If $x\in I$ then all the elements below $x$ are also in $I$.  Removing these $\al_s+1$ elements from $I$ yields an ideal $J$ in $F(\al_1,\ldots,\al_{s-2},\al_{s-1}-1)$ and this is a bijection.  This accounts for the second term in the recursion.
\eprf

In order to make use of this lemma, we will have to consider the indices where the coefficients of a polynomial achieve their maximum.
Given $f=\sum_k a_k q^k$ we define the {\em set of maxima indices} as
$$
\mi(f) = \{k \mid a_k = m\}
$$
where $m$ is the maximum value of a coefficient of $f$.  Note that if $f$ is unimodal then $\mi(f)$ will be an interval of integers.  The next result is easy to prove so the demonstration is omitted.
\ble
\label{mi}
Let $f,g$ be unimodal polynomials and suppose that $\mi(f)\cap\mi(g) \neq\emp$.  Then $f+g$ is unimodal and 
$\mi(f+g)=\mi(f)\cap \mi(g)$.\hqed
\ele

We can now prove a stronger version of Theorem~\ref{long} for the first segment.
\bth
\label{al_1}
If $\al=(\al_1,\al_2,\ldots,\al_s)$ satisfies
$$
\al_1\ge \al_2+\al_3+\cdots+\al_s
$$
then $r(q;\al)$  is unimodal with
$$
\mi(r(q;\al))=[\al_2+\al_3+\cdots+\al_s,\ \al_1].
$$
\eth
\bprf
 We will induct on $s$ and, for a given $s$, on $\al_2+\ldots+\al_s$.  We will only provide details for the case when $s$ is odd.
So we will apply the first recursion of Lemma~\ref{al_sRec}.  Using the induction hypothesis on the first term of the recursion, we see that 
$$
\mi(r(q;\be))=[\al_2+\al_3+\cdots+\al_s-1,\ \al_1]
$$
where $\be= (\al_1,\ldots,\al_{s-1},\al_s-1)$.  Similarly, for the second term we have
$$
\mi(q^{\al_s+1}\cdot r(q;\ga))= [\al_2+\al_3+\cdots+\al_s,\ \al_1+\al_s+1]
$$
where $\ga=(\al_1,\ldots,\al_{s-2},\al_{s-1}-1)$.  Now
$$
[\al_2+\al_3+\cdots+\al_s-1,\ \al_1] \cap  [\al_2+\al_3+\cdots+\al_s,\ \al_1+\al_s+1] = [\al_2+\al_3+\cdots+\al_s,\ \al_1]\neq\emp
$$
by our hypothesis on $\al$.  So Lemma~\ref{mi} applies and we are done.
\eprf

We have the same result for the last segment.
\bth
\label{al_s}
If $\al=(\al_1,\al_2,\ldots,\al_s)$ satisfies
$$
\al_s\ge \al_1+\al_2+\cdots+\al_{s-1}
$$
then $r(q;\al)$  is unimodal.
\eth
\bprf
If $s$ is even, then $L(\al)\iso L(\al^r)$.  So the result follows from the previous theorem. 
If $s$ is odd then, by Lemma~\ref{al^r}, $L(\al)\iso L(\al^r)^*$.  Since taking the dual just reverses the coefficients of the rank polynomial, we are again done by the previous result.
\eprf

\section{Arbitrary long segments}
\label{lms}

We will now complete the proof of Theorem~\ref{long}.  We will also provide an inductive way of creating posets which satisfy Conjecture~\ref{centered}.  We begin by locating the ranks of maximum size in the lattice coming from a poset with a long segment.

\ble
\label{MaxRk}
Let $\al=(\al_1,\ldots,\al_s)$ and $n=\#F(\al) $.  Suppose that for some $t$ we have
\beq
\label{al_t}
\al_t> \sum_{i\neq t} \al_i.
\eeq
Then the maximum size of a rank of $L=L(\al)$ is $\ell=\#L(F')$ where $F'$ is the poset obtained by removing the elements of segment $t$ from $F=F(\al)$. And this maximum occurs at ranks $m+1$ through $n-m-1$ where $m=\#F'$.
\ele
\bprf
Any $I\in L$ has the form
$I=J\cup K$ where $J\in L(F')$ and $K\in L(S)$ with $S$ being segment $t$ of $F$.
So if $\# I = k$ is fixed, the maximum number of possible ideals at this rank is the number of choices for $J$ since $S$ has only one ideal of each rank.  The maximum number of possible $J$ is $\ell$, which implies $r_k(\al)\le \ell$ for all $k$.  

We will now show that
$r_k(\al)=\ell$ for $m+1\le k\le n-m-1$ which will finish the proof.  Note that $[m+1,n-m-1]\neq\emp$ because of assumption~\ree{al_t}.
Every $I$ at such a rank $k$ must contain $\min S$ since, if not, then $I\sbe F'$ which forces $\#I\le m$.
Similarly, $I$ can not contain $\max S$ since if it did then $\#I\ge n-m$.  From these two statements, it follows that 
$I=J\cup K$ is an ideal for any choice of $J$, and for any $J$ there is a unique $K$ yielding an ideal at rank $k$ since $S$ is a chain.
So the number of choices for $I$ equals the number of choices for $J$ which is $\ell$, as desired.
\eprf

We will now prove Theorem~\ref{long}, noting that we have already demonstrated a slight strengthening of it in Theorems~\ref{al_1} and~\ref{al_s} for the first and last sements.
\bth
Suppose $\al=(\al_1,\ldots,\al_s)$ satisfies, for some $t$,
\beq
\label{al_tMid}
\al_t> \sum_{i\neq t} \al_i.
\eeq
Then $r(q;\al)$ is unimodal.
\eth
\bprf
We adopt the notation of the statement and proof of Lemma~\ref{MaxRk}.  According to this lemma, it suffices to
find a matching in $L$ from rank $k$ into rank $k+1$ for $k\le m$ and from rank $k$ into rank $k-1$ for $k\ge n-m$.

First consider the case $k\le m$.  We claim that if $\#I = k$ and $y$ is the coatom on segment $S$ then $y\not\in I$.  For suppose 
$y\in I$.  Then all the elements under $y$ on $S$ are in $I$ which implies $\#I\ge\al_t= n-m-1$.  However, 
$m=\# F' =\sum_{i\ne t} \al_i$ so equation~\ree{al_tMid} forces $n-m-1\ge m+1$ or $n\ge 2(m+1)$.
It follows that $\#I\ge 2(m+1)-m-1 = m+1$ which is a contradiction. 

Now given $I$ with $\#I\le m$, we match this ideal with $I\cup\{y_I\}$ where $y_I$ is the smallest element on $S$ not in $I$.  From the previous paragraph $y_I\le y$ so that $I\cup\{y_I\}$ is still an ideal.  And this is clearly a matching since $I$ determines $y_I$ uniquely.

The case $k\ge n-m$ is similar.  For this range of $k$ one can show that $x\in I$ where $x$ is the atom on $S$.  So one can match $I$ with $I-\{x_I\}$ where $x_I$ is the largest element on $S$ which is in $I$.  The details are similar to the first case and so left to the reader.
\eprf

We now give an inductive method for proving that Conjecture~\ref{centered} holds.  To do so, we must first investigate the finer structure of $L(\al)$ where $\al$ has a long segment.  For any ranked poset $P$, let 
$$
P_k =\{x\in P \mid \rk x \le k\}
$$
and
$$
P^k =\{x\in P \mid \rk x \ge k\}.
$$
\ble
\label{iso}
Let $\al=(\al_1,\ldots,\al_s)$, $F=F(\al)$, and $L=L(\al)$.  Also let $n=\#F$ and $m=\#F'$ where $F'$ is as in Lemma~\ref{MaxRk}.  Suppose that for some $t$ we have
$$
\al_t> \sum_{i\neq t} \al_i.
$$
Let $G=F(\be)$ and $M=L(\be)$ where
$$
\be=(\al_1,\ldots,\al_{t-1},\al_t+1,\al_{t+1},\ldots,\al_s).
$$
Then we have isomorphisms
$$
L_{n-m-1}\iso M_{n-m-1}\qmq{and} L^{n-m}\iso M^{n-m+1}.
$$
\ele
\bprf
We will construct an isomorphism $f:M_{n-m-1}\ra L_{n-m-1}$.  Note that $F$ can be obtained from $G$ by contracting (in the sense of graph theory) the edge of the Hasse diagram of $G$ containing the maximum vertex $z$ and the coatom $y$ it covers on the $t$th segment into a single maximum vertex.   
Doing this contraction on the second segment in Figure~\ref{F(2,3,1)} yields the fence in Figure~\ref{contract}.
But if $J\in M$ has $\rk J \le n-m-1$, then the inequality for $\al_t$ ensures that $J$ does not contain either of these two largest vertices.  So this induces a bijection  between such $J$ and the corresponding $I\in L_{n-m-1}$ which is order preserving in both directions, being essentially the identity map.

We can also define an isomorphism $g:M^{n-m+1}\ra L^{n-m}$.  This is similar to the construction of $f$.  
Let $w$ and $x$ be the minimum element and atom on $S$, respectively.  Note that if $J\in M^{n-m+1}$ then $w,x\in J$.  So contracting $x$ into $w$ yields the isomorphism.
\eprf

\bfi
\begin{tikzpicture}
\fill(0,0) circle(.1);
\fill(1,1) circle(.1);
\fill(2,2) circle(.1);
\fill(3,1) circle(.1);
\fill(4,0) circle(.1);
\fill(5,1) circle(.1);
\draw (0,0)--(2,2)--(4,0)--(5,1);
\draw (0,-.5) node{$x_1$};
\draw (1,.5) node{$x_2$};
\draw (2,1.5) node{$z$};
\draw (3,.5) node{$x_5$};
\draw (4,-.5) node{$x_6$};
\draw (5,.5) node{$x_7$};
\end{tikzpicture}
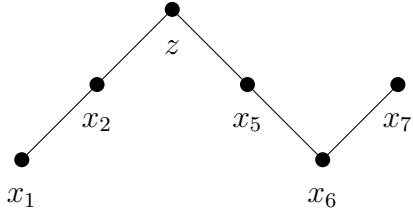
\capt{The fence obtained by contracting $x_3$ and $x_4$ in Figure~\ref{F(2,3,1)} \label{contract}}
\efi

We now have everything in place to state our inductive criterion for checking whether a poset with a long segment satsifies Conjecture~\ref{centered}.
\bth
Assume the hypotheses and notation of Lemma~\ref{iso}.  If $L$ has an NCD  then so does $M$.
Furthermore, if the NCD of $L$ is symmetric, top centered, or bottom centered then the NCD of $M$ has the same property.
\eth
\bprf
Suppose we have such an NCD of $L$.  We will keep the notation of the proof of Lemma~\ref{iso}.  Then, by this lemma, this NCD lifts to NCDs of $M_{n-m-1}$ and $M^{n-m+1}$.  What remains is to show how to connect these two NCDs through rank $n-m$ of $M$ to form an NCD of the whole poset.  
Let $R$, $R'$, and $R''$ denote ranks $n-m-1$ , $n-m$, and $n-m+1$ of $M$, respectively.  Note that by Lemma~\ref{MaxRk} we have $\#R=\#R'=\ell$.  As noted in the previous proof, no element of $R$ contains the coatom of the long
sement, $S$.  So there is a perfect matching from $J\in R$ to $J'\in R'$ by letting $J'=J\cup\{a\}$ where
where $a=\min(S-J)$.  To find the chains which continue to $R''$, consider $I=f(J)$ to
 see if it is connected in the NCD of $L$ to some $I''$ in the next rank up.  If so, then $J'$ will be covered by $J''=g^{-1}(I'')$  and so the sequence $J,J',J''$ will serve to connect the corresponding chains in $M$.  
It is easy to see that this construction will preserve whether the NCD is symmetric, top centered, or bottom centered which
completes the proof.
\eprf

By  induction on the length of the long segment, we immediately get the following result.  Note that showing one lattice has an NCD of a certain form immediately gives an infinite family of lattices with NCDs of the same form.
\bco
Let $\al=(\al_1,\ldots,\al_s)$ where
\beq
\label{plus1}
\al_t= 1 + \sum_{i\neq t} \al_i
\eeq
for some $t$.  If $L=L(\al)$ has an NCD  then so does $M=L(\be)$ where
$$
\be=(\al_1,\ldots,\al_{t-1},\al_t+a,\al_{t+1},\ldots,\al_s)
$$
for any $a\ge 0$.  Furthermore, if the NCD of $L$ is symmetric, top centered, or bottom centered then the NCD of $M$ has the same property. \hqed
\eco

Now one can use a computer to verify the following.
\bco
Let $F=F(\al_1,\ldots,\al_s)$ where
$$
\al_t> \sum_{i\neq t} \al_i
$$
where
$$
\sum_{i\neq t} \al_i \le 5.
$$
Then $L=L(\al)$ satisfies Conjecture~\ref{centered}.\hqed
\eco

Note that with a little more care, one can replace condition~\ree{plus1} with
$$
\al_t=\sum_{i\neq t} \al_i
$$
when $t=1$ or $t=s$.

\section{At most three segments}
\label{amts}

This section is devoted to a proof of Theorem~\ref{<=3}.  The result is trivial if the fence $F(\al)$ has one segment since in that case the poset is a chain and so its lattice is as well.  As we will see in the next result, the case of two segments is also easy.  But for three, we will have to use a modified version of the famous Greene-Kleitman symmetric chain decomposition of the Boolean algebra of all subsets of a finite set~\cite{GK:svs}.

\bth
Let $\al=(a,b)$.  Then $L(\al)$ is bottom interlacing.
\eth
\bprf
Let $z$ be the maximum element of $F=F(\al)$.  Also let $C,D$ be the first and second segments of $F$, respectively, with $z$ removed.
Note that any ideal of $F$ except $F$ itself is the disjoint union of an order ideal of $C$ and an order ideal of $D$.  It follows that
$$
L(\al)\iso (L(C)\times L(D))\oplus \{\oh\}
$$
where $\times$ is poset product and $\oplus$ is ordinal sum which makes every element of the product less than $\oh$.
Since $C$ and $D$ are chains, so are $L(C)$ and $L(D)$.  And there is a well-known symmetric decomposition $C_1, C_2,C_3,\ldots,$ of a product of two chains.  So in $L(\al)$ we have $\cen C_i = (n-1)/2$ for all $i$ where $n=\#F$.  Let $C_1$ be the chain containing the minimum element of $L(C)\times L(D)$.  Then, by symmetry, $C_1$ must also contain the maximum element of the product.  So $C_1'= C_1\cup\{\oh\}$ is a saturated chain in $L(\al)$ with $\cen C_1' = n/2$.  It follows that $C_1', C_2, C_3, \ldots$ is the desired bottom-interlacing NCD.
\eprf

To finish the proof of Theorem~\ref{<=3}, we will use the idea of a Greene-Kleitman core.  To define this object, let 
$w=w_1 w_2\ldots w_n$ be a sequence (or word) of zeros and ones.  The {\em Greene-Kleitman (GK) core of $w$}, $\GK(w)$, 
is a set of pairs of indices formed as follows.  If $w_i=0$ and $w_{i+1}=1$ then $(i,i+1)\in\GK(w)$.  We continue to add pairs $(i,j)$ to 
$\GK(w)$ as long as $w_i=0$ and $w_j=1$, where $i<j$ and all the elements between $w_i$ and $w_j$ are already in pairs of the core.  For example, if
$$
w= 110001011000111
$$
then
$$
\GK(w) = * * *\ \widehat{0\widehat{01}\widehat{01}1}* \widehat{0\widehat{01}1}\ *
$$
where elements not in the GK core have been replaced by stars, and pairs in the core are indicated by the hats.  Writing out the pairs themselves gives
$$
\GK(w) = \{(4,9),\ (5,6),\ (7,8),\ (11,14),\ (12,13)\}. 
$$
We will refer to the elements in the pairs of $\GK(w)$ as {\em matched}.  

To apply this idea to fences, we will have to modify the GK core.  And to do that it will be convenient to think of a fence $F$ as a partial order on $[n]=\{1,2,\ldots,n\}$ where, as usual, $n=\#F$.  When doing this, it will be important to distinguish $i\le j$ which is the usual  total order on the integers and $i\lte j$ which will be an order relation in $F$.  So consider $F=F(a,b,c)$ as the fence with covering relations
$$
b+c+1 \lt b+c+2 \lt \ldots \lt a+b+c+1 \gt b \gt b-1 \gt \ldots \gt 1 \lt b+1 \lt b+2 \lt \ldots \lt b+c.
$$
In other words label the second segment except for its maximum element with the elements of the interval $[1,b]$ from bottom to top.  Then label the elements of the third segment (except its minimum which has already been labeled) bottom to top with $[b+1,b+c]$.  Finally label the complete first segment with $[b+c+1,a+b+c+1]$ again from bottom to top.  Note that this labeling is a linear extension of $F$.  This labeling for the fence $F(2,3,1)$ is showing in Figure~\ref{labeling}.  

\bfi
\begin{tikzpicture}
\fill(0,1) circle(.1);
\fill(1,2) circle(.1);
\fill(2,3) circle(.1);
\fill(3,2) circle(.1);
\fill(4,1) circle(.1);
\fill(5,0) circle(.1);
\fill(6,1) circle(.1);
\draw (0,1)--(2,3)--(5,0)--(6,1);
\draw (0,.5) node{$5$};
\draw (1,1.5) node{$6$};
\draw (2,2.5) node{$7$};
\draw (3,1.5) node{$3$};
\draw (4,.5) node{$2$};
\draw (5,-.5) node{$1$};
\draw (6,.5) node{$4$};
\end{tikzpicture}
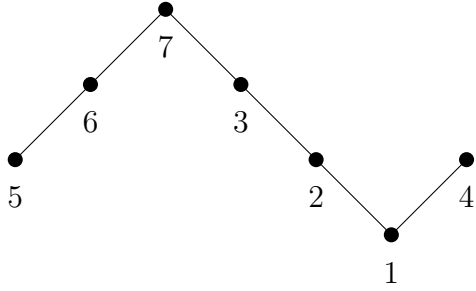
\capt{The labeling of $F(2,3,1)$ \label{labeling}}
\efi

Now associate with any subset $S\sbe F$ a word $w=w_S=w_1\ldots w_n$ where $w_i$ is one or zero depending on whether $i\in S$ or $i\not\in S$, respectively, and $n=\#F$.  Suppose $I\sbe F$ is an ideal and $w=w_I$.  
Since specifying $I$ and specifying $w$ are equivalent, we will often go back and forth without mention.
We will need a lemma about where $n$ lies with respect to the GK core of $w$ for ideals of certain three-segment fences
\ble
\label{unmatched}
Let $F=F(a,b,c)$ where $a\ge c$, $n=\#F$, and $I$ be an ideal of $F$.  Then $n$ is unmatched.
\ele
\bprf
Suppose, towards a contradiction, that $n$ is matched.  Since $n$ is the largest index in $w$, it must be in a pair $(i,n)$ and so $n\in I$.  Since $I$ is an ideal, this forces $[1,b]\cup[b+c+1,a+b+c+1]\sbe I$.  But then every element of $[b+c+1,a+b+c+1]$ must be matched with some element of $[b+1,b+c]$.  This contradicts the fact that $a\ge c$.
\eprf

Given an ideal, $I$, call an element $f\in F$ {\em frozen} if it is unmatched in $\GK(w)$ and 
there exists $(i,j)\in GK(w)$ with $f \gt i$ or $f\lt j$.  Note that since $(i,j)$ is in the GK core we must have $i\not\in I$ and $j\in I$.  So 
$f\gt i$ implies $f\not\in I$, since $I$ is an ideal, and $w_f=0$.  Similarly, if $f\lt j$ then $f\in I$ and $w_f=1$.  
Whether an element is frozen or not depends upon the ideal under consideration, but we will make sure $I$ is clear from context.
Finally, define the 
{\em core} of $w=w_I$ by
$$
\core w = \GK(w) \cup \{f\in F \mid \text{$f$ is frozen}\}.
$$
We say that elements of $F$ not in $\core w$ are {\em free}.
For example, suppose $F=F(2,3,1)$ and $I=\{1,4,5\}$.  Then we first compute the GK core as indicated by the hats in
$$
w= 1 \widehat{0 \widehat{0 1} 1} 0 0.
$$ 
Since $(3,4)\in \GK(w)$ and $1\lt 4$ we have that $1$ is frozen in $w$.
Similarly, $(3,4)\in\GK(w)$ and $7\gt 3$ implies that $7$ is frozen.  One can also check that $6$ is free so that
$$
\core w = 1 \widehat{0 \widehat{0 1} 1} * 0
$$
or
$$
\core w = \{1,\ (2,5),\ (3,4),\ 7\}.
$$
We will need some facts about frozen elements.
\ble
\label{frozen}
Let $F=F(a,b,c)$ where $a\ge c$, $n=\#F$, and $I$ be an ideal of $F$.
\ben
\item[(a)]  If $f\in F$ is frozen then $f=1$ or $f=n$.
\item[(b)]  If $1$ is frozen then $w_1=1$.  If $n$ is frozen then $w_n=0$.
\item[(c)]  If $1$ is frozen, then so is $n$.
\een
\ele
\bprf
Let the segments corresponding to $a,b,c$ be $A,B,C$, respectively.
Suppose $(i,j)\in\core w=w_I$.  So $i\not\in I$, $j\in I$ and $i<j$.  It follows that $i,j$ must be in different segments since the labeling of $F$ is a linear extension and $I$ is an ideal.  The order of labeling the segments and $i<j$ then imply that there are only three different possibilities: $i\in B$ and $j\in C$, or $i\in B$ and $j\in A$, or $i\in C$ and $j\in A$.

We now consider what elements are frozen for each of the three positions of $i,j$.  Suppose first that $i\in B$ and $j\in C$.
Note that we can not have $i=n$ since $i<j$.  
If we have an element $f$ satisfying $f\gt i$ then for $f$ to be frozen it must be unmatched.  But all the elements of $B-\{n\}$ above $i$ come between $w_i$ and $w_j$ and so are in pairs of $\GK(w)$.  So the only element of $B$ which could be frozen is $n$.
In fact,  $n$ must be frozen since $n\gt i$ forces $w_n=0$ and a final zero is never matched.  In a similar way, we see that if $f\lt j$ then it can only be frozen if $f=1$ and that this element must indeed be frozen.

In a similar manner, one shows that if $f$ is frozen when  $i\in B$ and $j\in A$ then $f=n$ and that $n$ is always frozen in this case.
Finally, we consider  $i\in C$ and $j\in A$.  As before, we can show that no $f\gt i$ is frozen.  As far as the $f\lt j$, note that $j\neq n$ because of the previous lemma.  So all such $j$ are on segment $A$ and now one can proceed as before to show that no $f\lt j$ are frozen.

Summarizing the three cases in order, if $(i,j)\in\core w$ then it freezes: $1$ and $n$, or just $n$, or no element of $F$.  And in all three cases if $1$ or $n$ is frozen then $w_1=1$ and $w_n=0$.  This completes the proof of the lemma.
\eprf

We are now prepared to prove the last case of Theorem~\ref{<=3} which we restate here for convenience.
\bth
\label{3seg}
If $\al=(a,b,c)$ then $L(\al)$ admits a CD which is symmetric, top heavy, or bottom heavy depending on whether $a=c$, $a>c$, or $a<c$, respectively.
\eth
\bprf
We will consider that case $a\ge c$.  This is without loss of generality since if $a<c$ then $\al^r$ has its first part larger than its third.

Let $I$ be an ideal of $F=F(\al)$ and let $w=w_I$.  The free elements of $I$ must form a subsequence of $w$ where all the ones precede all the zeros since, by the previous lemma, if this were not true then there would be a pair of free elements which could be matched.   Let $w'$ be the sequence formed from $w$ by replacing the left-most free zero $w_k$ (if any)  with a one.  We will show that the corresponding $I'\sbe F$ is an ideal and that $\core w'=\core w$.

Suppose, towards a contradiction, that $I'$ is not an ideal .  Then we must have $i\lt k$ where $i\not\in I$.  Now $i$ can not be free in $w$ because of the choice of $k$.  Also $i$ can not be paired since, if it were, then $k$ would be frozen and not free.  Finally, $i$ can not be frozen since Lemma~\ref{frozen} (a) would imply that $i=1$, which contradicts Lemma~\ref{frozen} (b).  So $i$ does not exist and $I'$ is an ideal.

Now suppose, towards another contradiction, that  $\core w'\neq\core w$.  Then it must be that $k$ is an element of $\core w'$.  But $k$ can not be paired in $w'$ since, by the choice of $k$, all $w_i=0$ with $i<k$ are already in pairs of $\core w\sbe\core w'$.  So $k$ must be frozen in $w'$.  By the previous lemma and the fact that $w_k'=1$ we must have $k=1$.  However $w_1=0$ and $I$ is an ideal, so $I$ contains no elements on either the second or third segment of $F$.  This contradicts the fact that, from the proof of Lemma~\ref{frozen}, the only way $w_1'=1$ can be frozen is if there is an $(i,j)\in\core w'$ where $j\gt 1$ is on the third segment.

Consider what happens if we form a sequence $w''$ from $w$ by turning the right-most free one (if any) into a  zero.  We can similarly prove that the resulting $I''$ is an ideal with the same core as $w$.

Now given $\ka=\core w_I$ where $I$ is an ideal, we form a corresponding chain in $L=L(\al)$, $C_\ka$, as 
follows.  The minimum element on $C_\ka$ is  the word $v$ which has core $\ka$ and all free elements equal to zero.  Now turn the free $0$s in $v$ into $1$s one at a time, moving left to right until they are all $1$s.  From what we have proved, the elements of $C_\ka$ are exactly the ideals of $L$ with core $\ka$.  And the chain is saturated by construction.  From this discussion it is clear that the $C_\ka$ as 
$\ka$ runs over all possible cores form a CD of $L$.  Furthermore, from the previous lemma, every chain either has $\cen C_\ka = n/2$ if $1$ and $n$ are either both frozen or both not, or has $\cen C_\ka = (n-1)/2$ if $n$ is frozen but $1$ is not.  Thus this is a bottom-heavy CD.

Finally, we consider what happens when $a=c$.  In that case $F\iso F^*$ so that $L\iso L^*$.  It follows that $r(q;\al)$ is symmetric and so the bottom-heavy CD we have constructed must, in fact, be symmetric.
\eprf

\section{$d$-divided posets}
\label{ddf}

We will now give a simplified proof of a theorem of Gansner~\cite{gan:loi} showing that the posets he considered, which we will call $d$-divided, have top-interlacing NCDs.  We note that Gansner did not state his result in terms of the top-interlacing concept since it had not been defined when his paper was written.  But it is easy to derive this property from his proof.

Let $n,d$ be positive integers.  Divide $n$ by $d$ to obtain 
\beq
\label{div}
n=qd+r
\eeq
where $0\le r <d$.  We then define the corresponding 
{\em $d$-divided partially ordered set}, $P_{n,d}$, to be the poset on $[n]$ with covering relations
$$
d\gt d-1 \gt \ldots\gt 1\lt 2d\gt 2d-1 \gt \ldots \gt d+1 \lt 3d \gt \ldots \lt n\gt n-1\gt\ldots\gt qd+1.
$$
So $P_{n,d}$ has segments with $d$ elements alternating with segments of $2$ elements, ending with a segment of $r$ elements.
And the segments of length $d$ and $r$ are labeled left to right from the bottom to the top of each segment.  So, again, the labeling is a linear extension.  For example, $P_{10,4}$ is displayed in Figure~\ref{P_104}.

\bfi
\begin{tikzpicture}
\fill(0,0) circle(.1);
\fill(0,1) circle(.1);
\fill(0,2) circle(.1);
\fill(0,3) circle(.1);
\fill(2,0) circle(.1);
\fill(2,1) circle(.1);
\fill(2,2) circle(.1);
\fill(2,3) circle(.1);
\fill(4,2) circle(.1);
\fill(4,3) circle(.1);
\draw (0,3)--(0,0)--(2,3)--(2,0)--(4,3)--(4,2);
\draw(-.4,0) node{$1$};
\draw(-.4,1) node{$2$};
\draw(-.4,2) node{$3$};
\draw(-.4,3) node{$4$};
\draw(1.6,0) node{$5$};
\draw(1.6,1) node{$6$};
\draw(1.6,2) node{$7$};
\draw(1.6,3) node{$8$};
\draw(3.6,2) node{$9$};
\draw(3.6,3) node{$10$};
\end{tikzpicture}
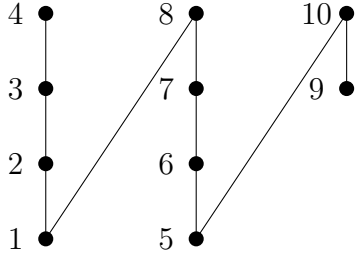
\capt{The $3$-divided poset $P_{10,4}$ \label{P_104}}
\efi

Now given any ideal $I$ of $P_{n,d}$, we have a corresponding word $w=w_I$ just as in the previous section.  So we can form its Green-Kleitman core $\GK(w)$ as before.  To define what the core is for this poset, we need to consider $p\in P_{n,d}$ which is the lowest element on the right-most chain of length $d$.  In our example, $p=5$.  Consider the subword $v=w_p w_{p+1} \ldots w_n$ of $w$.
We now define the core of $w$ to be
$$
\core w =\case{\GK(w)}{if $n$ is unmatched in $\GK(v)$,}{\GKh(w)}{else,}
$$
where $\GKh(w)$ is a modification of $\GK(w)$ defined as follows.  If $n$ is in a pair of $\GK(w)$, then we must have $w_n=1$ since this is the last element of the word.  This forces $w_p=1$ because $p\lt n$ and $I$ is an ideal.  We now consider $w_p$ as paired with itself and form the rest of $\GKh(w)$ around it in the usual way.  Returning to our running example, suppose $I=\{1,2,5,6,9,10\}$ so that
$$
w = 1100110011.
$$
Then 
$$
v= 110011
$$
with
$$
\GK(v) = **\widehat{0\ \widehat{0\ 1}\ 1}.
$$
Since $w_{10}$ is matched in $\GK(v)$, we form $\GKh(w)$ by considering $w_5$ as matched with itself which will be indicated by it having a hat of its own.  So
$$
\core w=\GKh(w) = ***\ \widehat{0\ \widehat{1}\ 1}\ \widehat{0\ \widehat{0\ 1}\ 1}.
$$
As in the previous section, we call the elements of $[n]$ {\em matched} or {\em free} with respect to $w$ depending on whether they are in $\core w$ or not, respectively.  We note that if $r=0$ in~\ree{div} then $n$ is always unmatched in $v$ and so in this case 
$\core w=\GK(w)$.

We are now ready to prove that $L(P_{n,d})$ always has an NCD which is either top-heavy or symmetric.  One can translate these into instances where Conjecture~\ref{centered} is true by noting that the dual $P_{n,d}^*$ is a fence and $L(P_{n,d}^*)\iso L(P_{n,d})^*$.
Our proof will parallel that of Theorem~\ref{3seg}.
\bth[\cite{gan:loi}]
The lattice $L(P_{n,d})$ has an NCD which is symmetric or top heavy depending on whether $d$ divides $n$ or not, respectively.
\eth
\bprf
Let $I$ be an ideal of $P_{n,d}$ and let $w=w_I$.  As in the proof of Theorem~\ref{3seg}, the free elements of $w$ form a subsequence of ones followed by a subsequence of zeros.  We will consider what happens in passing to $w'$ which is obtained from $w$ by replacing the rightmost free one $w_i$ (if any) by a zero.  We will show that $w'$ is the word of an ideal $I'$ and that 
$\core w'=\core w$.

We suppose $I'$ is not an ideal and derive a contradiction.  So there must be $j\in I$ with $j\gt i$.  It follows that $j$ and $i$ are on the same segment $S$.  We will now consider three cases depending on the location of $S$.

First suppose that $S$ is one of the segments with $d$ or $r$ elements as defined by~\ree{div}.  Take $j$ minimal with respect to these restrictions.  It follows that $j-1\not\in I$.  If $i<j-1$ then we immediately have a contradiction since $j\in I$, $j-1\not\in I$, 
$j-1\lt j$, and $I$ is an ideal.  So $i=j-1$ which means $w_{i+1}=w_j=1$.  But $w_i$ is the rightmost free one, so $w_{i+1}$ must be matched.  This contradicts the fact that a matched one can only have matched ones between it and its paired zero.

We now consider what happens when $S$ is one of the segments of length one, first supposing that it is not the last such segment.  
So all of the element on the segment $S'$ of length $d$ containing $j$ must be in $I$.  But there are only $d-1$ other elements between $w_i$ and $w_j$.  Thus at least one of the elements of $S'$ must be paired with an element before $i$, contradicting the fact that $i$ is free in $w$.

Finally, we look at the case when $S$ is the last segment of length one.  So $j=n$ and $i=p$ in the notation of the discussion before this theorem.  If $n$ is not matched in $\GK(v)$ then $\core w =\GK(w)$ where $n$ is matched since $p$ was the rightmost free one.  Now using an argument as in the previous paragraph, we get a contradiction.  But if $n$ is matched in $\GK(v)$, we get an immediate contradiction because $p$ is in $\core w$ and so not free.

 Now that we know $I'$ is an ideal, we wish to show $\core w'=\core w$.  But the only way $\core w'$ could change is if the new zero in position $i$ becomes paired with a one.  And  this one would have to be free and  in a position to the right of $i$, contradicting the fact that $i$ was the rightmost free one.

Similarly one can show that changing the leftmost free zero to a one in $I$ results in an ideal with the same core.  The rest of the proof is almost exactly like the last two paragraphs of the demonstration of Theorem~\ref{3seg}.  The only difference is that now a chain $C_\ka$ will have 
$\cen C_\ka = n/2$ or $(n+1)/2$ where the latter case happens whenever there is a one paired with itself in $\ka$.
\eprf
 
\section{Two more approaches}
\label{tma}

We end with two other possible methods for proving the conjectures in this paper.  While we have not been able to make them bear fruit, we hope that someone else will be more successful.

\medskip

{\bf 1.  An explicit formula.}  There is an explicit formula for $r(q;\al)$ in terms of powers of $q$ and $q$-integers 
$[n]_q=1+q+\cdots + q^{n-1}$ where $n\ge0$.  This expression appeared in a more complicated form in Claussen's thesis~\cite{cla:epp} and formed part of the basis for his proof that Conjecture~\ref{MGO} is true for any fence with at most four segments.  It is easy to prove by toggling in a manner similar to that of the demonstration of Lemma~\ref{al_sRec}.  But now one toggles on subsets of 
\beq
\label{Z(al)}
Z(\al)=\{z_1,z_2,\ldots,z_u\}
\eeq
which is the set  of  maxima of $F(\al)$ written in order from left to right.  We will state the theorem without proof since the demonstration does not contain any new ideas.  A similar result holds for an even number of segments.
\bth[\cite{cla:epp}]
Let $\al=(\al_1,\ldots,\al_{2u-1})$ with $Z(\al)$ given by~\ree{Z(al)}.  Then
$$
r(q;\al) = \sum_{Z\sbe Z(\al)}  q^{\#Z} r(Z,1) r(Z,2) \cdots r(Z,u)
$$
where 
$$
r(Z,1) =
\begin{cases}
q^{\al_1} & \mbox{if $z_1\in Z$,}\\
[\al_1+1]_q & \mbox{if $z_1\not\in Z$,}
\end{cases}
$$
and
$$
r(Z,i)=
\begin{cases}
q^{\al_{2i-2}+\al_{2i-1}-1} & \mbox{if $z_{i-1},z_i\in Z$,}\\
q^{\al_{2i-2}}[\al_{2i-1}]_q & \mbox{if $z_{i-1}\in Z$ and $z_i\not\in Z$,}\\
q^{\al_{2i-1}}[\al_{2i-2}]_q & \mbox{if $z_{i-1}\not\in Z$ and $z_i\in Z$,}\\
1+ q[\al_{2i-2}]_q[\al_{2i-1}]_q &\mbox{if $z_{i-1},z_i\not\in Z$.}
\end{cases}
$$
for $i>1$.\hqed
\eth

One can now write out $r(q;\al)$ for any desired $\al$ and try to use inductive arguments on these expressions.  It is easy to prove that
multiplying by $[n]_q$ preserves the interlacing conditions.  Unfortunately, addition is another matter and this is where the difficulty lies.

\medskip

{\bf 2.  Lexicographic CDs.}  There are two hurdles to using the Greene-Kleitman approach to constructing NCDs.  The first is finding a suitable labeling of the poset.  The other is coming up with an appropriate modification of $\GK(w)$ when $w$ is taken with respect to the labeling.  Here we propose another method which, given any labeling of any finite labeled poset $P$, produces a chain decomposition of the lattice of ideals $L=L(P)$ which may or may not be nested.

We construct the chains $C_1, C_2, C_3,\ldots$ of the CD  as follows.  Suppose $C_1,\ldots,C_{i-1}$ have been constructed. 
Since $P=[n]$ as sets, we can consider any ideal $I$ of $P$ as a subset of $\{1,\ldots,n\}$ and we will not make any distinction between an ideal and its subset.  So given two ideals, we can compare them in the lexicographic order on subsets.  Now we form $C_i$ by starting with the unique ideal $I_0$ which has minimum rank and is also lexicographically least among all elements of $L'=L-(C_1\cup\cdots\cup C_{i-1})$.
We now consider all ideals of $L'$ which cover $I_0$ and take the lexicographically least of them to be the next element $I_1$ on $C_i$.  We continue in this manner until we come to an ideal which has no cover in $L'$ at which point $C_i$ terminates.  We have the following conjecture which we have verified  for all fences $F$ with $\#F\le 7$.
\bcon
For every $\al$ there is a labeling of $F(\al)$ such that the correponding lexicographic CD is an NCD.
\econ

There are a couple of stumbling blocks we have encountered in trying to prove this conjecture.  One is guessing what the proper labeling of $F(\al)$ should be.  It seems that for a given $\al$ there are many labelings which work.  So narrowing the field down is where the difficulty lies.  The second problem is trying to prove that you have a CD for a certain family of fences once you have identified a suitable labeling.  Induction is the obvious proof technique, but a small change in the fence can make the lexicographic CD change in what seems to be difficult ways. 

Still, we are optimistic that much progress can be made on these conjectures.  We invite the reader to work on them.



\end{document}